\documentclass[11pt]{article}

\usepackage{amsfonts}
\usepackage{amssymb}
\usepackage{amsmath}
\usepackage{amsthm}
\usepackage{graphicx}
\usepackage{float}
\usepackage{bm}
\usepackage{color}
\usepackage{subcaption}
\usepackage{cite}
\usepackage{url}
\usepackage[ruled,linesnumbered]{algorithm2e}
\usepackage{mathtools}

\theoremstyle{plain}

\theoremstyle{remark}

\title{Evaluation of resonances: adaptivity and AAA rational
  approximation\\ of randomly scalarized boundary integral resolvents}

\author{Oscar P. Bruno\footnote{Computing and Mathematical Sciences, California Institute of Technology, Pasadena, CA, 91125 USA, obruno@caltech.edu, msantana@caltech.edu)}
  \and Manuel A. Santana\footnotemark[1]
\and Lloyd N. Trefethen\footnote{School of Engineering and Applied Sciences, Harvard University, Cambridge, MA 02138 USA 
  (trefethen@seas.harvard.edu)}
}

\usepackage{amsopn}

\DeclareMathOperator*{\argmin}{arg\,min}

\def \C{\mathbb{C}}

\def \inv{^{-1}}
\usepackage{todonotes}
\usepackage[T1]{fontenc}
\usepackage{listings}
\date{}
\begin{document}

\maketitle

\begin{abstract}
  This paper presents a novel algorithm, based on use of rational
  approximants of a randomly scalarized boundary integral resolvent in
  conjunction with an adaptive search strategy and an exponentially
  convergent secant-method termination stage, for the evaluation of
  acoustic and electromagnetic resonances in open and closed cavities;
  for simplicity we restrict treatment to cavities in two-dimensional
  space. The desired cavity resonances (also known as ``eigenvalues''
  for interior problems, and ``scattering poles'' or ``complex
  eigenvalues'' for exterior and open-cavity problems) are obtained as
  the poles of associated rational approximants; both the approximants
  and their poles are obtained by means of the recently introduced AAA
  rational-approximation algorithm. In fact, the proposed
  resonance-search method applies to any nonlinear eigenvalue problem
  associated with a given function $F: U \to \C^{d\times d}$, wherein,
  denoting $F(k) = F_k$, a complex value $k$ is sought for which
  $F_kw = 0$ for some nonzero $w\in \C^d$. For the scattering problems
  considered in this paper, which include interior, exterior and open
  cavity problems, $F_k$ is taken to equal a spectrally discretized
  version of a Green function-based boundary integral operator at
  spatial frequency $k$. In all cases, the scalarized resolvent is
  given by an expression of the form $u^* F_k\inv v$, where
  $u,v \in \C^d$ are fixed random vectors. The proposed adaptive
  search strategy relies on use of a rectangular subdivision of the
  resonance search domain which is locally refined to ensure that all
  resonances in the domain are captured. The approach works equally
  well in the case in which the search domain is a one-dimensional
  set, such as, e.g., an interval of the real line, in which case the
  rectangles used degenerate into subintervals of the search domain.
  A variety of numerical results are presented, including comparisons
  with well-known methods based on complex contour integration, and a
  discussion of the asymptotics that result as open cavities approach
  closed cavities---in all, demonstrating the accuracy provided by the
  method, for low- and high-frequency states alike.
\end{abstract}
\section{Introduction\label{intro}}

We are concerned with the problem of evaluation of resonances
supported by open and closed cavities and other scattering structures,
which are obtained as solution pairs $(u,k)$ of the eigenvalue problem
\begin{equation}
  \label{eq:helm} \Delta u + k^2u = 0
\end{equation}
with eigenfunction $u$ and eigenvalue $-k^2$, posed on an interior or
exterior domain $\Omega$, with homogeneous boundary conditions of
e.g.\ Dirichlet, Neumann, or other types.  Once approximated by
discretized versions of the problem's boundary integral operators
(which is done in this paper on the basis of the open- and
closed-curve integral equation algorithms
of~\cite{bruno2012second,lintner2015generalized,COLTON:2019}, see
also~\cite{bruno2013high}), the resonance-search problem is reduced to
the solution of a related Nonlinear Eigenvalue Problem (NEP) for a
certain function $F: U\to \C^{d\times d}$, wherein, denoting
$F(k) = F_k$, a complex value $k$ is sought for which
\begin{equation}\label{NEP}
  F_kw = 0\quad \mbox{for some nonzero}\quad w\in \C^d.
\end{equation}
This contribution focuses on scattering problems, for which the
function $F$ in~\eqref{NEP} provides a discrete approximation of the
associated boundary-integral operator.  However the method
is general and, indeed, other NEPs unrelated to boundary integral
operators are considered in this paper as well.

The proposed approach seeks the desired resonant values $k$, for
which~\eqref{NEP} holds, as poles of the randomly scalarized resolvent
\begin{equation}
  \label{eq:scalarized}
  S(k) = u^* F_k\inv v, \quad\mbox{where}\quad u,v \in \C^d\quad \mbox{are fixed random  vectors}.
\end{equation}
The desired resolvent poles are obtained as poles of rational
approximants of~\eqref{eq:scalarized}; both the rational approximants
and their poles are produced numerically by means of the recently
introduced AAA rational-approximation algorithm~\cite{AAA}.  The
proposed eigensolver additionally incorporates an adaptive search
strategy and a secant-method termination stage. The adaptive search
strategy relies on use of a rectangular subdivision of the resonance
search domain which is locally refined to ensure that all resonances
in the domain are captured. The adaptivity approach works equally well
in the case in which the search domain is a one-dimensional set, such
as, e.g., an interval of the real line, in which case the rectangles
used degenerate into subintervals of the search domain. The
secant-method termination stage, in turn, is an important element in
the proposed algorithm, which enables (i)~Exponentially fast
convergence to near machine precision accuracy starting from AAA-based
results of lower accuracy; (ii)~Reliable error estimation; and,
(iii)~A capability to reliably screen the spurious eigenvalues that
can (rarely) be produced by the AAA algorithm.  In all, the overall
proposed approach is simple, easy to implement and rapidly convergent,
and it requires limited computation besides the embarrassingly
parallelizable evaluation of the scalarized resolvent at various
wavenumbers $k$. A variety of numerical results presented in this
paper demonstrate the character of the proposed approach: the method
yields highly accurate approximations of scattering resonances and
solutions of other NEPs, even in cases involving high frequencies.

A significant literature has developed in recent years in connection
with the solution of NEPs. As discussed in the survey
article~\cite{guttel2017nonlinear}, solution methods include root
finding methods, contour integration methods, and methods based on
linearization of $F_k$, all of which have been applied to the
computation of resonances \cite{EldarBoundary, alves2024wave,
  el2020rational,misawa2017boundary,steinbach2012convergence,
  zhao2015robust}. In turn, a set of methods for the NEP that, like
the present paper, rely on use of AAA rational approximation, have
recently been developed~\cite{lietaert2022automatic,guttlescalar1,
  guttelscalar2}, and specifically, the
contributions~\cite{guttlescalar1, guttelscalar2} apply the AAA
algorithm to the scalarized resolvent~\eqref{eq:scalarized}. But in
these contributions the AAA algorithm is applied in a manner different
to the one we use: in those contexts a rational approximation is
employed to produce a linearization of $F_k$ whose eigenvalues
approximate the desired eigenvalues, whereas the present paper
directly uses the poles of the rational approximant of $S$ as
approximants of the desired eigenvalues. Closer to our work is the
AAA-based algorithm introduced in \cite{betz2024efficient}, which
considers the transmission of plane waves by a periodic dieletric
system. In that approach, numerical solutions of the transmission
problem are obtained by means of the finite element method and from
such solutions the rational approximant of the coefficient of
transmission across the structure is produced---whose poles near the
real axis are indicative of the resonant character of the structure,
in that they incorporate information concerning some of the
structure's complex eigenvalues.

In view of the strengths of the various algorithmic components it
incorporates, the proposed algorithm is flexible, accurate and
efficient. The algorithm's use of rational approximations allows for
utilization of scalarized-resolvent data points on arbitrary sets of
frequencies, enabling, in particular, the use of arbitrarily
distributed data points on the boundaries of rectangles---which
greatly facilitates the design of the adaptive strategy proposed in
this paper. Further, use of simple intervals of the real line suffices
for evaluation of real eigenvalues. The secant-method termination
stage, finally, provides significant benefits concerning accuracy and
reliability, as mentioned in points (i)-(iii) above.

This paper is organized as follows. Sections~\ref{eigs_int_ops}
and~\ref{num_int} review the integral-equation and associated
numerical schemes used in this paper to represent solutions of the
Helmholtz equation~\eqref{eq:helm} for open and closed two-dimensional
domains.  Section \ref{aaa} provides a brief description of the AAA
algorithm. The overall proposed approach for the solution of NEPs is
then described in Section~\ref{sec:NEP}. A variety of numerical
results presented in Section~\ref{num_ex} include NEPs unrelated to
Laplace eigenvalues, comparisons with well-known methods based on
complex contour integration, illustrations concerning low- and
high-frequency Laplace eigenvalue problems for open and closed
cavities, as well as the asymptotics that result as open cavities
approach closed cavities---in all, demonstrating the accuracy provided
by the method even for low- and high-frequency states alike.

\section{Eigenvalue problems, Green Functions and Integral Operators\label{eigs_int_ops}}
We consider eigenvalue problems of the form~\eqref{eq:helm}, posed in
open two dimensional spatial domains $\Omega$ with smooth boundaries
$\Gamma$, and with homogeneous boundary conditions on $\Gamma$.  Three
types of spatial domains are considered in this paper, namely, domains
$\Omega$ equal to: (a)~The complement $\Gamma^c$ of an open arc
$\Gamma$ in $\mathbb{R}^2$; (b)~The region interior to a closed curve
$\Gamma$ in $\mathbb{R}^2$; and, (c)~The region exterior to a closed
curve $\Gamma$ in $\mathbb{R}^2$.  For definiteness this paper mostly
concerns eigenvalue problems under homogeneous Dirichlet
boundary conditions
\begin{equation}
  \label{eq:dirichlet}
  u|_\Gamma = 0
\end{equation}
for each of these domain-types. Homogeneous Neumann and Zaremba
boundary conditions can be handled
similarly~\cite{lintner2015generalized,EldarBoundary,akhmetgaliyev2017regularized};
in fact, results of Neumann problems produced by the proposed
algorithm are briefly mentioned in Section~\ref{sec:gapsize}.

Our treatment of the problem~\eqref{eq:helm}--\eqref{eq:dirichlet} is
based on use of the two dimensional Helmholtz Green function
$G_k(x,y)\coloneq \frac{i}{4} H^1_0(k|x-y|)$ (where $H^1_0$ denotes the
Hankel function of the first kind of order $0$) and the associated
single-layer potential representation
\begin{equation}\label{ansatz} 
  u(x)= \int_{\Gamma} G_k (x,y) \psi(y) \, ds_y, \quad x \in \Omega,
\end{equation}
of the eigenfunction $u$ in terms of a surface density $\psi$. In view
of the well known~\cite{COLTON:1983,mclean2000strongly} continuity of
the single layer potential $u = u(x)$ as a function of
$x\in\mathbb{R}^2$, up to and including $\Gamma$, we consider the
boundary integral operator $F_k: H^{-1/2}(\Gamma) \to H^{1/2}(\Gamma)$
(resp.\  $F_k: \widetilde{H}^{-1/2}(\Gamma) \to H^{1/2}(\Gamma)$) given
by the expression
\begin{equation}
\label{modelB}
F_k[ \psi](x)  =  \int _{\Gamma}G_k(x,y)\psi (y)ds_y, \quad x \in \Gamma
\end{equation}
on a closed (resp.\ open) smooth curve $\Gamma$.
See~\cite{mclean2000strongly} and~\cite{lintner2015generalized} and
references therein for detailed definitions of the Sobolev spaces
$H^{\pm 1/2}$ and $\widetilde{H}^{1/2}$ relevant to the closed and open-arc single
layer operators, respectively.

As suggested above, the $k$-dependent operator~\eqref{modelB} can be
used to tackle the interior, exterior and open-arc eigenvalue problems
under consideration. Indeed, for any smooth open arc or closed curve
$\Gamma$ and for any given density function $\psi$ defined on
$\Gamma$, we have~\cite{mclean2000strongly} (i)~The function $u$ given
by the representation formula~\eqref{ansatz} is a solution of
equation~\eqref{eq:helm} for all $x\in\Gamma^c$; and, (ii)~The
function $u\ne 0$ satisfies the Dirichlet boundary
condition~\eqref{eq:dirichlet} if and only if $\psi\ne 0$ is a
solution of the equation $F_k[\psi] = 0$. It can accordingly be shown
that the resolvent operator $(F_k)^{-1}$ is an analytic function of
$k$ for $\Im k >0$, and that that given a complex number $k = \mu$
with $\Im \mu\leq 0$, the number $-\mu^2$ is an eigenvalue of the
problem~\eqref{eq:helm}--\eqref{eq:dirichlet} if and only if the value
$k = \mu$ is a pole of the resolvent operator $(F_k)^{-1}$ as a
function of $k$. This fact is established
in~\cite[Prop.\ 7.10]{taylorbook} for closed-curve interior problems
(for which $\mu$ in the lower half-plane $\Im \mu \le 0$ must actually
be real) and for closed-curve exterior problems (for which $\mu$ must satisfy $\Im \mu < 0$).  As indicated in
what follows, the corresponding result for open-arc problems can be
established along lines similar to those
in~\cite[Prop.\ 7.10]{taylorbook}.

A critical element in the extension of these results to the open-arc
case is the injectivity of the mapping $\psi \to u$, which, according
to equation~\eqref{ansatz}, maps functions defined on $\Gamma$ to
functions defined in $\mathbb{R}^2$. The corresponding closed-curve
injectivity result for the exterior problem is established
in~\cite{taylorbook} by showing that if $\psi$ satisfies the equation
$F_k(\psi)=0$, then the associated function $u$ given
by~\eqref{ansatz} is a Laplace eigenfunction in the interior of
$\Gamma$ and that, therefore, the corresponding eigenvalue $-\mu^2$
must be real, and by subsequently making use of the jump relations for
the single- and double-layer potentials across $\Gamma$. The
corresponding closed-curve result for the interior problem is
established similarly. For the open-arc problem we have no equivalent
of the interior region, but the injectivity result can be established
nevertheless, simply by using the jump relation for the normal
derivative of the single-layer potential. The equivalence between
Laplace eigenvalues $-\mu^2$ with $\mu$ in the lower half-plane and
poles $k = \mu$ of the resolvent $(F_k)^{-1}$ then follows for the
open arc case in a manner analogous to that put forth
in~\cite[Sec.\ 7]{taylorbook} by relying on the second-kind formulation
for open problems introduced
in~\cite{bruno2012second,lintner2015generalized}.

In sum, noting that, without loss of generality, the search for values
of $k$ satisfying the eigenvalue
problem~\eqref{eq:helm}--\eqref{eq:dirichlet} may be restricted to the
lower half-plane $\Im k \leq 0$, the eigenvalues may be sought as real
poles $k = \mu$ of the resolvent operator $(F_k)^{-1}$ for the
eigenvalue problem in the interior of a closed curve $\Gamma$, and as
complex poles $k = \mu$ of the same operator, with $\Im \mu <0$, for
either the eigenvalue problem in the exterior of a closed curve
$\Gamma$ or for the complement of an open arc $\Gamma$. The associated
eigenfunctions $u$ then result via equation~\eqref{ansatz} with
density (or, for multiple eigenvalues, densities) $\psi$ in the
nullspace of the operator $F_k$. In other words, the integral equation
setting just described reduces the eigenvalue
problem~\eqref{eq:helm}--\eqref{eq:dirichlet} for interior and
exterior closed-curve and open arc problems to an NEP for the
single-layer operator~\eqref{modelB} with values of $k$ restricted to
the lower half-plane.

\section{Numerical Instantiation of the Integral Operator
  \boldmath${F_k}$\label{num_int}}
The numerical implementations utilized in this paper for the Green
function-based integral operator~\eqref{modelB} are based on the
discretization methods presented in~\cite[Sec.\ 3.6]{COLTON:2019}
(resp.\ \cite[Secs.\ 3.2, 5.1]{bruno2012second}) for closed (resp.\
open) curves in the plane. For simplicity, our computational examples
restrict attention to smooth curves $\Gamma$ and boundary conditions
of Dirichlet type, although related methods are
available~\cite{akhmetgaliyev2017regularized,bruno2012second,COLTON:2019}
that enable corresponding treatments for non-smooth
boundaries~\cite{akhmetgaliyev2017regularized,COLTON:2019} as well as
Neumann, Robin and Zaremba boundary
conditions~\cite{akhmetgaliyev2017regularized,bruno2012second}. As
indicated in~\cite{bruno2012second}, in particular, the numerical
implementation of the integral operator~\eqref{modelB} for open arcs
$\Gamma$ requires consideration of the edge singularities that are
incurred by the solutions $\psi$ of problems of the form
$F_k[\psi] = f$ even for functions $f$ defined on $\Gamma$ which do
not contain such singularities.

For both open- and closed-curve problems the
methods~\cite{COLTON:2019,bruno2012second} discretize the single layer
operator $F_k$ on the basis of Nystr\"om-type
methodologies---utilizing a sequence of points along the curve
$\Gamma$ for both integration and
enforcement of the equation. For closed curves $\Gamma$ the
discretization is produced as the image under the curve
parametrization of a uniform grid in the parameter interval
$[0,2\pi]$; in the case of open curves the discretization results as
the image of a Chebyshev mesh in the parameter interval $[-1,1]$. In
both cases the unknown functions $\psi$ are expressed in terms of
Fourier-based expansions in the parameter intervals, which are then
integrated termwise by reducing each integrated term to evaluation of
explicit integrals. The edge singularities of the function $\psi$ in
the open-arc case are tackled by explicitly factoring out the singular
term: using a smooth parametrization $\mathbf{r} = \mathbf{r}(t)$ of
the open curve $\Gamma$ ($-1\leq t\leq 1$) and writing
$\psi(\mathbf{r}(t)) =\phi(\mathbf{r}(t))/\sqrt{1-t^2}$, it follows
that $\phi$ is a smooth function. Upon introduction of a cosine change
of integration variables two desirable effects occur, namely, the
function $\phi$ is converted into a periodic and even function which
may be expanded in a cosine series with high accuracy, and, at
the same time, the square root term in the denominator is cancelled by
the Jacobian of the change of variables. The method is completed by
exploiting explicit expressions for the integral of products of a
logarithmic kernel and the cosine Fourier basis functions. As
illustrated in~\cite{bruno2012second} and other contributions
mentioned above, these methodologies can produce scattering solutions
with accuracies near machine precision on the basis of relatively
coarse discretizations, even for configurations involving high spatial
frequencies.

In particular, these procedures produce highly accurate numerical
approximations of the integral operator $F_k$
in~\eqref{ansatz}---which we exploit in the context of this paper to
produce accurate numerical evaluations of eigenvalues and
eigenfunctions. As indicated in Section~\ref{sec:NEP}, the poles of (a
randomly scalarized version of) this integral operator, which, per the
discussion in Section~\ref{eigs_int_ops}, correspond to Laplace
eigenvalues in the various cases considered, are then obtained as
poles of associated AAA rational approximants. The
corresponding eigenfunctions are obtained via consideration of a
Gaussian elimination-based de-singularization procedure described in
Section~\ref{sec:NEP}.  For added accuracy, the methods in that
section propose two alternatives, namely, the use of either iterated
AAA rational approximants on one hand, or application of the secant
method after an initial evaluation of poles via the AAA approach, on
the other hand. In practice we have observed that, without exceptions,
accuracies near machine precision are obtained for both eigenvalues
and eigenfunctions on the basis of the overall proposed methodology; a
few related illustrations are presented in Section~\ref{num_ex}.

\section{Rational Approximants and the AAA Algorithm\label{aaa}}

The AAA algorithm is a greedy procedure for the construction of a
rational approximant to a given complex-valued function $f$ on the
basis of its values on an $N$-point set $Z_N\subset\mathbb{C}$ (full
details can be found in~\cite{AAA}).  Given the set
$\{(z,f(z))\ |\ z\in Z_N\}$, the algorithm proceeds by selecting a
sequence of points $z_j\in Z_N$, starting with some point $z_1$, which
in principle can selected arbitrarily, but which the Matlab
implementation~\cite{chebfun} takes as a $z\in Z_N$ for which the
function value $f(z)$ is farthest from the mean of the set
$\{ f(z)\ |\ z\in Z_N\}$. The remaining points are then selected
inductively. Once points $z_j\in Z_N$ ($1\leq j \leq m$) have been
chosen, with corresponding function values $f_j \equiv f(z_j)$, for a
suitably chosen vector $w^m = (w^m_1,\dots,w^m_m)$ of complex weights
$w^m_j$ satisfying $\sum |w_j|^2 = 1$, the procedure to obtain
$z_{m+1}$ starts by constructing the barycentric-form rational
function
\begin{equation}\label{rat_fun}
  r(z) = \frac{n^m(z)}{d^m(z)} = \sum_{j=1}^m \frac{w^m_j f_j}{z - z_j} \bigg/   \sum_{j=1}^{m} \frac{w^m_j}{z - z_j},
\end{equation}
where, as suggested by the notation used, $n^m(z)$ and $d^m(z)$ denote the
numerator and denominator in the right-hand expression
in~\eqref{rat_fun}. The vector $w^m$ is selected as follows: calling
$A_m = \{(\widetilde{w}_1,\dots,\widetilde{w}_m)\in \mathbb{C}^m\ | \sum |\widetilde{w}_j|^2 = 1
\}$, $w^m$ is defined as a minimizer of the least-squares problem
\[
  w^m = \argmin_{\widetilde w\in A_m} \sum_{z\in
    Z^m_N}|f(z)d_{\widetilde w}(z) -n_{\widetilde w}(z)|^2,
\]
where $Z^m_N = Z_N \setminus \{z_j\ |\ j=1,\dots, m\}$,
$n_{\widetilde w}(z) = \sum_{j=1}^m \frac{\widetilde w_j f_j}{z -
  z_j}$ and
$d_{\widetilde w}(z)= \sum_{j=1}^m \frac{\widetilde w_j}{z - z_j}$.
Once the minimizer $w^m$ has been computed, $z_{m+1} $ is defined as
a point $z \in Z^m_N$ for which $|f(z) -n^m(z)/d^m(z)|$ is
maximum relative to $\max\{ f(z)\ |\ z\in Z_N\}$. The algorithm terminates when this maximum is less than or
equal to a specified tolerance, and the last rational
function~\eqref{rat_fun} obtained as part of the $z_j$ selection
process provides the desired rational approximant.

All of the numerical illustrations presented in this paper
utilize the AAA implementation included with Chebfun~\cite{chebfun}.

\section{Solution of NEPs\label{sec:NEP}}
The discussion in Sections~\ref{eigs_int_ops} and~\ref{num_int}
reduces the eigenvalue problem \eqref{eq:helm}--\eqref{eq:dirichlet}
to NEPs for the single-layer operator~\eqref{modelB} (and the
corresponding discrete approximate operator introduced in
Section~\ref{num_int}), with values of $k$ restricted to the lower
half-plane. This section presents numerical algorithms for the
solution of this NEP and, indeed, of general NEPs for which the
resolvent operator $(F_k)^{-1}$ is a meromorphic function of $k$.

As indicated in Section~\ref{intro}, the proposed NEP algorithm
obtains the eigenvalues $k$ as the poles of a AAA rational
approximant associated with the randomly scalarized resolvent $S(k)$
in equation~\eqref{eq:scalarized}.  While, in principle, values of the
scalarized resolvent $S(k)$ on any given subset of the complex plane
may be used to obtain eigenvalues, throughout this paper we restrict
attention to algorithms based on use of values $S(k)$ for $k$ on a
given curve $\mathcal{C}$ in the complex plane. Both open and closed
curves $\mathcal{C}$ may be used, such as e.g. the closed curves equal
to the boundaries of either a rectangle or a circle, or the open curve
equal to an interval $[a,b]$ contained in the real axis. Proceeding on
the basis of such data, Algorithm~\ref{alg:one} below evaluates either
eigenvalues contained on the union $\widetilde{\mathcal{C}}$ of
$\mathcal{C}$ and its interior, if $\mathcal{C}$ is a closed curve, or
eigenvalues along the curve $\mathcal{C}$, if $\mathcal{C}$ is an open
curve. Algorithm~\ref{alg:two}, in turn, seeks to find all such
eigenvalues, namely, all eigenvalues contained within
$\widetilde{\mathcal{C}}$ if
$\mathcal{C}$ is a closed curve, and all eigenvalues along
$\mathcal{C}$ if $\mathcal{C}$ is an open curve.

\vspace{0.3cm}
\begin{algorithm}[H]
  Select a set $Z_N= \{k_1,\dots, k_N\}$ which is a subset of an open curve $\mathcal{C}$ or a closed curve and its interior $\mathcal{\Tilde{C}}$\\
  Choose random vectors $u,v \in \C^{d}$\\
  \For{$j = 1, \dots ,N$}{ $s_j = u^* F(k_j)\inv v $ }
  Compute a rational approximant $r(z)$ associated with the set $\{(k_j,s_j)\ |\ s_j= S(k_j), j=1,\dots,N\}$ using the AAA algorithm\\
  Return the poles of $r(z)$ in $\mathcal{C}$ or $\mathcal{\Tilde{C}}$
  \caption{Basic Algorithm}\label{alg:one}
\end{algorithm}
\vspace{0.3cm}

In detail, the pseudo-code Algorithm~\ref{alg:one} proceeds by first
evaluating the scalarized resolvent $S = S(k)$ at a suitably selected
set $Z_N= \{k_1,\dots, k_N\}$ of points along $\mathcal{C}$ to obtain
the set $G_N =\{(k_j,s_j)\ |\ s_j= S(k_j),\ j=1,\dots,N\}$. Using this
set of pairs Algorithm~\ref{alg:one} then obtains a rational
approximant $r = r(k)$ by means of the AAA algorithm, and, for a
closed curve $\mathcal{C}$, the poles of $r(k)$ contained in
$\widetilde{\mathcal{C}}$ (resp. for an open curve $\mathcal{C}$, the
poles on or sufficiently close to $\mathcal{C}$) are returned as
approximations of eigenvalues of $F_k$.

\vspace{0.3cm}
\begin{algorithm}[H]
    \SetKwFunction{EigRecursive}{eigadaptive}
    \SetKwProg{Fn}{Function}{}{}
    \DontPrintSemicolon 
    \KwIn{A set \textbf{C} which is either a rectangular region $\widetilde{\mathcal{C}}$ or real interval $\mathcal{C}$ and eigenvalues $e$ from Algorithm \ref{alg:one} applied to \textbf{C}}
    \Fn{\EigRecursive{$e$, \textbf{C}}}{
        Partition \textbf{C} dyadically in each dimension into sets $C_i$\\
        \For{each $C_i$}{
            Compute $n_i$, the number of eigenvalues $e$ in $C_i$\\
            Compute eigenvalues $\widetilde{e}$ from Algorithm \ref{alg:one} applied to $C_i$\\
            Compute $\widetilde{n}_i$, the number of eigenvalues $\widetilde{e}$ inside $C_i$\\
            \uIf{$n_i = \widetilde{n}_i$}{
            \Return $\widetilde{e}$
            }
            \uElse{
            \Return \EigRecursive{$\widetilde{e}$, $C_i$}
            }
        }
    }
\caption{Adaptive Algorithm} \label{alg:two}
\end{algorithm}
\vspace{0.3cm}

Typically Algorithm \ref{alg:one} finds all of the eigenvalues within
$\mathcal{C}$ or $\widetilde{\mathcal{C}}$, as relevant, provided the
number of such eigenvalues is not too large.  This observation suggests the development of an adaptive
version of Algorithm~\ref{alg:one}, which, as detailed in the
pseudo-code titled Algorithm~\ref{alg:two}, is specialized to searches
for eigenvalues within either (i)~A rectangular region
$\widetilde{\mathcal{C}}$, or, (ii)~A real interval
$\mathcal{C} = [a,b]$. The algorithm proceeds by finding all
eigenvalues within $\mathcal{C}$ or $\widetilde{\mathcal{C}}$, and
then reapplying Algorithm~\ref{alg:one} upon recursive dyadic
subdivision along each dimension. The recursion ends when the number
of poles obtained does not change upon subdivision. Clearly, the fact
that Algorithm~\ref{alg:one} is applicable to arbitrary curves is
highly beneficial in this context, as the use of rectangular regions
lends itself for the adaptive subdivision strategy used in
Algorithm~\ref{alg:two}.

It is important to note that the AAA algorithm may fail to accurately
produce the eigenvalues within $\mathcal{C}$ or $\widetilde{\mathcal{C}}$ if the set $G_N$
does not adequately represent $S(k)$ along the curve $\mathcal{C}$. An
inadequate representation generally manifests itself through the
appearance of false poles (see~\cite{AAA}) near $\mathcal{C}$, which
in practice may be detected via a convergence analysis as points are
added to the set $G_N$ and $N$ grows accordingly. Without exception,
in practice we have found that, provided the curve $\mathcal{C}$
encloses a sufficiently small number of eigenvalues, once convergence
to a fixed set of poles has occurred, the poles found within $\mathcal{C}$ or
$\widetilde{\mathcal{C}}$ correspond in a one-to-one fashion to all the
eigenvalues in that region. Starting from such poles, eigenvalues with
near machine precision accuracy can be obtained by means of one of two
possible methods, namely (i)~Use of a subsequent ``localized'' AAA
approximation applied to a few points around each eigenvalue obtained
in the initial application of AAA; or, (ii)~Use of the secant method
applied to $1/S(k) = 1/(u^*F_k^{-1}v)$ starting at each one of the
eigenvalues obtained in the initial application of AAA. In regard to
method~(i), we have found that use of localized AAA approximations
over a circle of radius $\rho_\mathrm{conv}=10^{-5}$ works well in
many cases, but generally tuning of the parameter $\rho_\mathrm{conv}$
is necessary for convergence and, indeed, it is difficult in practice
to determine whether convergence to a given tolerance has occurred,
except in test cases for which the eigenvalues are known
beforehand. Fortunately, method~(ii) does not present this difficulty
and, as it happens, it provides an additional significant
advantage. Indeed, in practice we have found that, without exception,
use of method~(ii) results in convergence to an eigenvalue, and,
clearly, the convergence history provides an error estimate for the
eigenvalue found---and we thus recommend this approach as a completion
procedure for each eigenvalue found.  The numerical examples in
Section~\ref{num_ex} illustrate the performance of
Algorithm~\ref{alg:two} augmented by means of each one of these
localized accuracy-improvement procedures.

Eigenvectors corresponding to each eigenvalue can finally be obtained
via evaluation, based on Gaussian elimination, of the nullspace of the
matrix $F_k$ at each eigenvalues $k$ obtained. In detail, using
Gaussian elimination with pivoting leads to an LU decomposition of the
form $F_k = L U$. Using this decomposition the nullspace of $F_k$ can
be computed by a simple two-step procedure consisting of (i)~Selection
of a set of canonical-basis vectors that are mapped to zero by the
rows in the matrix $U$ that are associated with the nonzero pivots;
and, (ii)~For each zero pivot in the matrix $U$, construction and
solution of the reduced systems that result from the elimination of
the rows and columns in $U$ containing the zero pivots, and with
right-hand sides equal to the negatives of each of the eliminated
columns but excluding the column element in the eliminated row. Null
vectors could also be computed by means of the singular value
decomposition for sufficiently small problems such as the two
dimensional examples considered in this paper, but the singular value
decomposition method would be problematic for larger problems such as
those arising from scattering in three dimensions, whereas, in view
of~\cite{reichel2005breakdown}, GMRES-based iterative methods related
to the LU-based approach mentioned above can be envisioned. An
alternative considered in~\cite[Sec.\ 2.2]{betz2024efficient} produces
eigenvectors on the basis of the rational approximation itself,
albeit, according to our experiments, at the expense of some loss of accuracy. 
In this contribution the aforementioned Gaussian elimination
procedure is used, as it requires no additional approximation after
application of the secant based refinement method, and as it results
in accuracies comparable to those enjoyed by the eigenvalues
themselves.

\section{Numerical Examples\label{num_ex}} A few numerical
illustrations of the proposed methodology are presented in what
follows, including a demonstration of the performance of Algorithm
\ref{alg:two} on NEPs unrelated to Laplace eigenvalues
(Section~\ref{sec:twoneps}); a set of examples concerning
low-frequency Laplace eigenvalues (Section~\ref{sec:lowfreq}) and a
comparison to complex contour integration methods
(Section~\ref{sec:contcompare}); an exploration of the rate at which
scattering poles associated with open arcs converge to the
corresponding interior eigenvalues as the opening closes
(Section~\ref{sec:gapsize}); and finally, a set of examples concerning
high-frequency Laplace eigenvalues and eigenfunctions
(Section~\ref{sec:highfreq}).

\subsection{Adaptivity and exhaustive eigenvalue
  evaluation}\label{sec:twoneps}
\begin{figure}[H]
    \centering
    \includegraphics[scale=0.4]{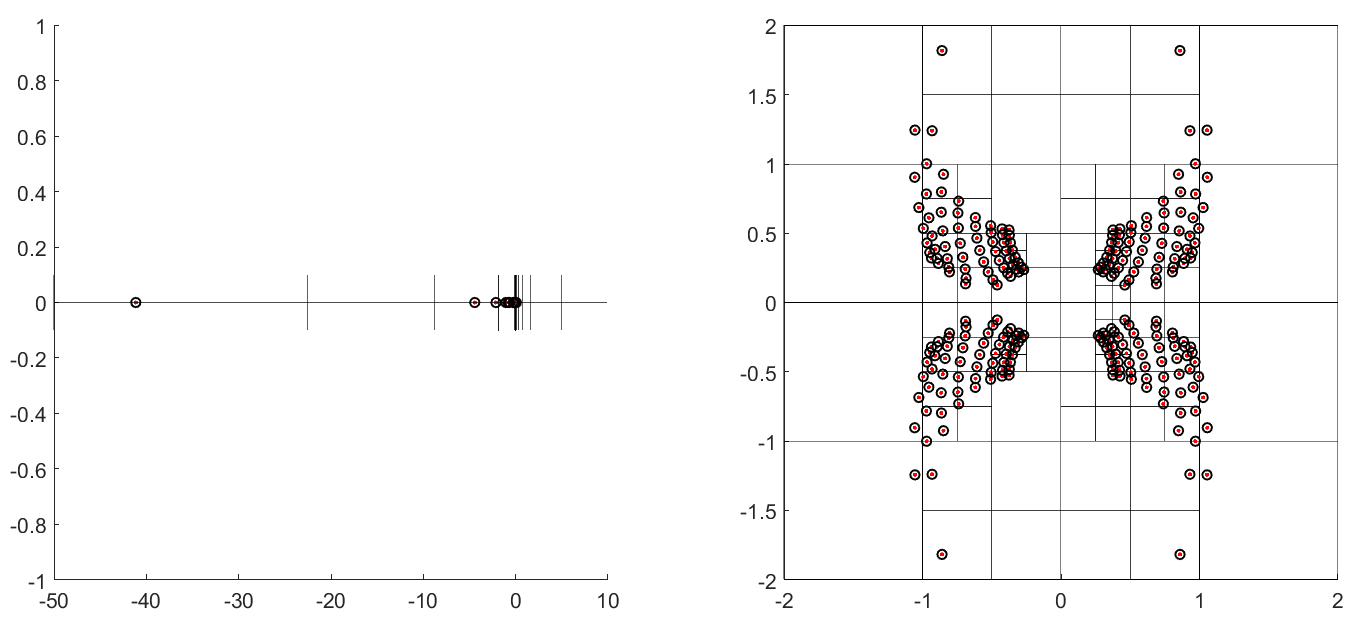}
    \caption{Demonstration of Algorithm \ref{alg:two} on two problems
      from the NLEVP set \cite{nlevp}. In both panels red points
      represent exact eigenvalues and black circles represent
      eigenvalues produced by the proposed numerical method. Black
      lines represent divisions related to the adaptive version of the
      algorithm. Left panel: The CD player problem, for which all $60$
      eigenvalues in the interval $[-50,5]$ were found to at least $7$
      digits. Right: The butterfly problem, for which all $256$
      eigenvalues were found to at least $10$ digits. Accuracy near
      machine precision was subsequently obtained by increasing the
      discretization in each subregion or by using methods refinement
      secant- or AAA-based methods described in Sections~\ref{sec:NEP}
      and~\ref{sec:lowfreq}}.
    \label{fig:NEP}
\end{figure}

This section illustrates the ability of Algorithm~\ref{alg:two} to
exhaustively evaluate the set of eigenvalues of a given NEP that are
contained in a given region in the complex plane (see
also the examples mentioned in the first paragraph of
Section~\ref{sec:highfreq}).  To do this, two well known problems
included in the NLEVP collection~\cite{nlevp} are considered, namely,
the ``CD player problem'' and the ``Butterfly problem''; the
corresponding eigenvalues are purely real in the first case and
complex but not real in the second case. The CD player problem is an
eigenvalue problem for a $60 \times 60$ matrix polynomial of the form
$F(k) = k^2I + kA_1 + A_0$ arising in the study of a CD-player control
task~\cite{cdplayer}. We restrict attention to the interval $[-50,5]$
on the real axis, within which the problem has $60$ eigenvalues with
absolute values as small as $2.23e-4$ and as large as $41.1399$; the
real-interval version of Algorithm 2 was used with a total number of
$300$ points $k_j$ along each subinterval.  The second problem is the
butterfly problem, in which eigenvalues of a $64\times 64$ matrix
polynomial $F(k) = \sum_{i=0}^4 k^iA_i$ are sought. The matrices $A_i$
are Kronecker products of linear combinations of the identity and
nilpotent Jordan blocks \cite{butterfly}. We tackle the butterfly
problem by employing Algorithm 2 applied initially to the box of side
length 4 centered at the origin. This problem has $256$ eigenvalues,
and $100$ values of $k$ were used on each side of the square. Figure
\ref{fig:NEP} displays the results produced by the proposed algorithm:
in each case the proposed adaptive algorithm obtains all of the
eigenvalues in the prescribed regions of the complex plane.

\begin{figure}
    \centering
    \includegraphics[scale=0.3]{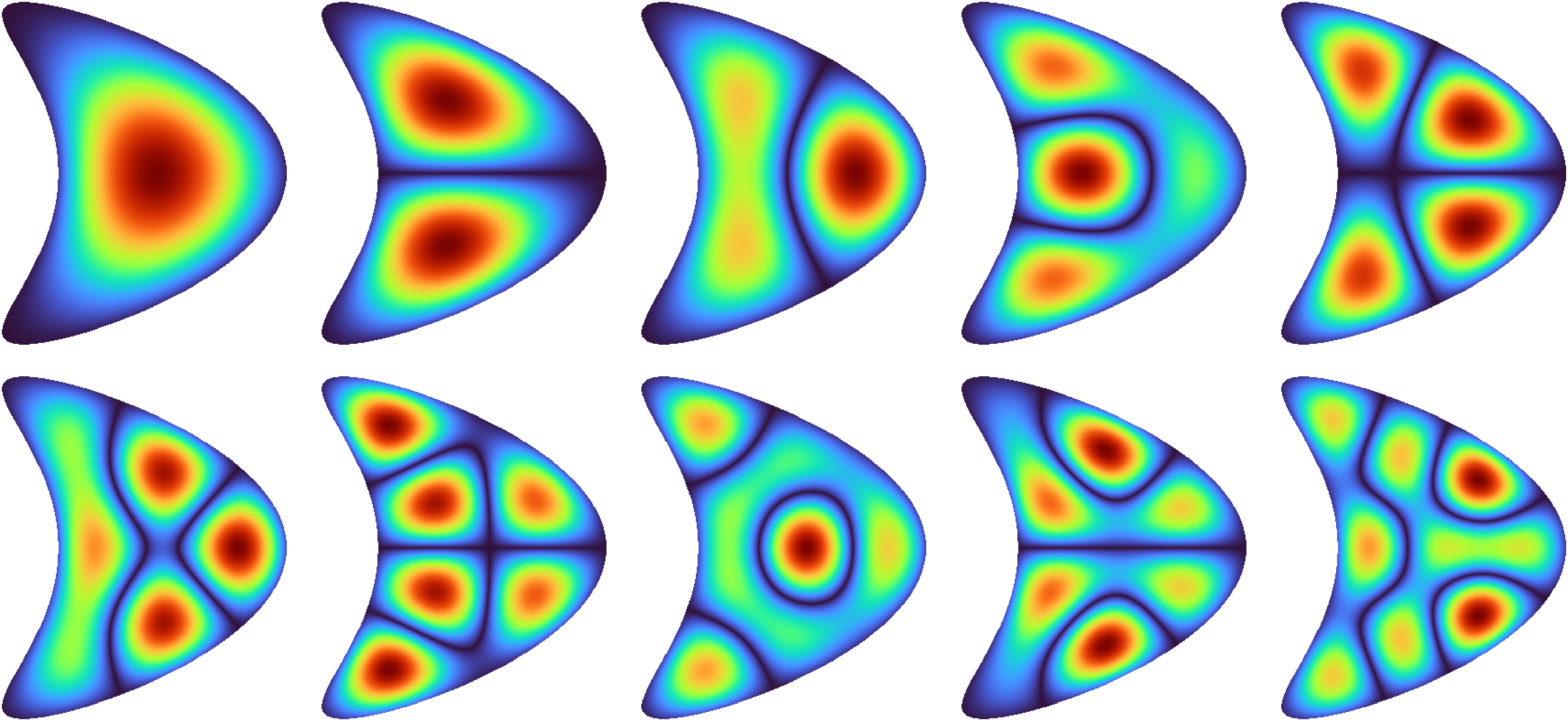}
    \caption{Low-frequency interior eigenvalue problems mentioned in
      Section \ref{sec:lowfreq}: the first ten interior eigenfunctions
      for the kite-shaped domain. The associated frequencies $k$ are
      listed in the left column of Table~\ref{tab:lowfreq_evals}. }
    \label{fig:kites}
\end{figure}

\subsection{Low-frequency eigenvalues}\label{sec:lowfreq}
\begin{figure}[h!]
    \centering
    \includegraphics[scale=0.38]{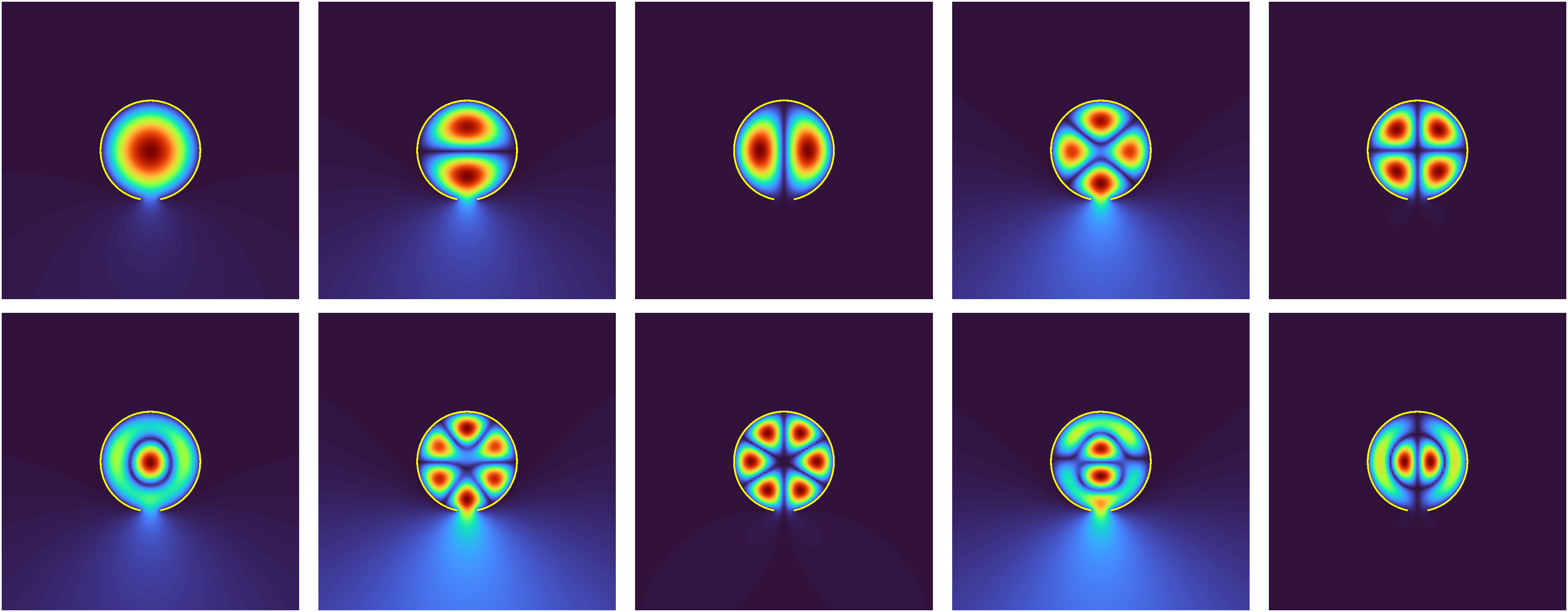}
    \caption{Low-frequency open arc eigenproblems discussed in
      Section~\ref{sec:lowfreq}: the first ten scattering poles of the
      open circle with imaginary part less than $-0.2i$ (ordered
      left-to-right and top-to-bottom with increasing values of the
      real part of the frequency). The associated frequencies $k$ are
      displayed in the right column of Table~\ref{tab:lowfreq_evals}.}
    \label{fig:circles}
  \end{figure}

This section illustrates the character of the proposed eigensolver in
the context of NEPs~\eqref{NEP} associated with low-frequency Laplace
eigenproblems~\eqref {eq:helm} with boundary
conditions~\eqref{eq:dirichlet}. Thus, Figures~\ref{fig:kites}
and~\ref{fig:circles} and Table~\ref{tab:lowfreq_evals} present
Laplace eigenfunctions and eigenvalues for two structures, namely, a
(closed) kite-shaped domain and a circular cavity with an aperture of
$\pi/8$ radians. All of these eigenvalues and eigenfunctions were
produced by means of Algorithm $2$ with an error (evaluated by
convergence studies) of the order of $\mathcal{O}(10^{-13})$. In the
first case the eigenvalues are real, and they are thus obtained by
means of Algorithm~$2$ with $\mathcal{C}$ equal to an interval in the
real axis. In the second case the eigenvalues lie near the real axis;
they were obtained using once again Algorithm~$2$, but, this time
using a rectangular curve $\mathcal{C}$ enclosing a finite section of
the strip between the real axis and the horizontal line
$\Im k = -0.2$; for such a thin strip the localization step was
applied by iteratively subdividing the domain along the real axis
only.

\subsection{Comparison with 
  Integral Algorithm~1 in~\cite{Beyn2012eigen} and the block SS method~\cite{asakura2009numerical}}\label{sec:contcompare}
It is most relevant in the context of this paper to compare the
proposed algorithm to other eigensolvers which, like ours, are based
on use of Green-function representations. Most often such algorithms
rely additionally on complex-contour integration strategies for
evaluation of eigenvalues~\cite{steinbach2012convergence,
  misawa2017boundary, leblanc2013solving, wieners2013boundary,
  hankelseperation, pommerenke2019computation,
  zolla2005foundations,bykov2012numerical, baum2005singularity}, and
we therefore select for our comparison the two best-known and
prototypical complex-contour eigensolvers, namely, the block SS
method~\cite{asakura2009numerical}, and ``Integral Algorithm~1''
in~\cite{Beyn2012eigen}. To avoid potential confusion with the
Algorithm~1 presented in Section~\ref{sec:NEP} above, in what follows
we refer to ``Integral Algorithm~1'' in~\cite{Beyn2012eigen} as
``Beyn~1''; for our comparisons we use the implementations of the
Beyn~1 and block SS algorithms provided in~\cite{guttel2017nonlinear}.

 \begin{table}[H]
    \centering
    \begin{tabular}{|c|c|}
    \hline
  Kite           & Open Circle\\
  \hline
  $2.209856180349$ &  $2.391850921204 - 0.000866833533i$\\
  $3.215653682128$ &  $3.785851440218 - 0.007551333804i$\\
  $3.528868275787$ &  $3.831519839558 - 0.000000810935i$\\
  $4.303831479675$ &  $5.066410135738 - 0.022753855105i$\\
  $4.371112240590$ &  $5.134599571714 - 0.000011845979i$\\
  $4.906513621606$ &  $5.486798760828 - 0.010839713761i$\\
  $5.291183742145$ &  $6.297659940294 - 0.044691641691i$\\
  $5.461743432329$ &  $6.377232306043 - 0.000071959651i$\\
  $5.736410337307$ &  $6.923647500434 - 0.056416692369i$\\
  $6.172352448525$ &  $7.015195622517 - 0.000013514954i$\\
  \hline
    \end{tabular}
    \caption{Computed frequencies $k$ (listed top-to-bottom in this
      table) such that $-k^2$ is an eigenvalue of the eigenvalue
      problem~\eqref{eq:helm} corresponding to the eigenfunctions
      displayed in Figures \ref{fig:kites} and \ref{fig:circles}
      (ordered left-to-right and top-to-bottom). All digits listed are
      believed correct on the basis of numerical convergence analyses.}
    \label{tab:lowfreq_evals}
  \end{table}

The block SS method relies on use of matrices containing $L$ columns
of random vectors by computing the product
\begin{equation}
  \label{eq:rand_matrix}
  U^* F_k\inv V\in \C^{L\times L}, \quad\mbox{where}\quad U,V \in \C^{d\times L}\quad \mbox{are matrices containing  $L$ random vectors as columns.}
\end{equation}
Additionally, per prescriptions in~\cite{asakura2009numerical}, a
number $2P$ of associated moments of the $k$ dependent $L\times L$
matrix in~\eqref{eq:rand_matrix} are computed and used to construct a
generalized eigenvalue problem whose solution provides numerical
values of the desired eigenvalues and eigenvectors. Here we
determine the number $2P$ of moments used (which according
to~\cite{guttel2017nonlinear} should be such that $PL$ is not smaller
than the number of eigenvalues contained in the contour, counting
multiplicities) on the basis of the SVD-based rank test introduced
in~\cite{Beyn2012eigen}. For simple eigenvalues, such as appear most
often for exterior scattering eigenvalues
problems~\cite{klopp1995generic}, it
suffices~\cite{asakura2009numerical} to choose $L = 1$ in the block SS
method. Note that with this selection of the parameter $L$, the block
SS algorithm is based on use of the scalarized
resolvent~\eqref{eq:scalarized}---exactly the same data utilized by
Algorithm 1. The Beyn~1 method, in turn, utilizes the product in
\eqref{eq:rand_matrix} with the matrix $U^*$ replaced by the
$d$-dimensional identity matrix. The number $L$ of columns used in the
Beyn~1 approach is selected on the basis of the aforementioned
SVD-based rank test; in particular, the value $L$ must at least equal
the number of eigenvalues within the contour (counting
multiplicities), and, to yield highly accurate results, it may need to
be increased by a small number, perhaps as small as one or two (in
accordance with the rank test), to account for the number of
eigenvalues outside the contour but sufficiently close to it.
\begin{figure}[h!]
        \hspace{-2cm} \includegraphics[scale=0.3]{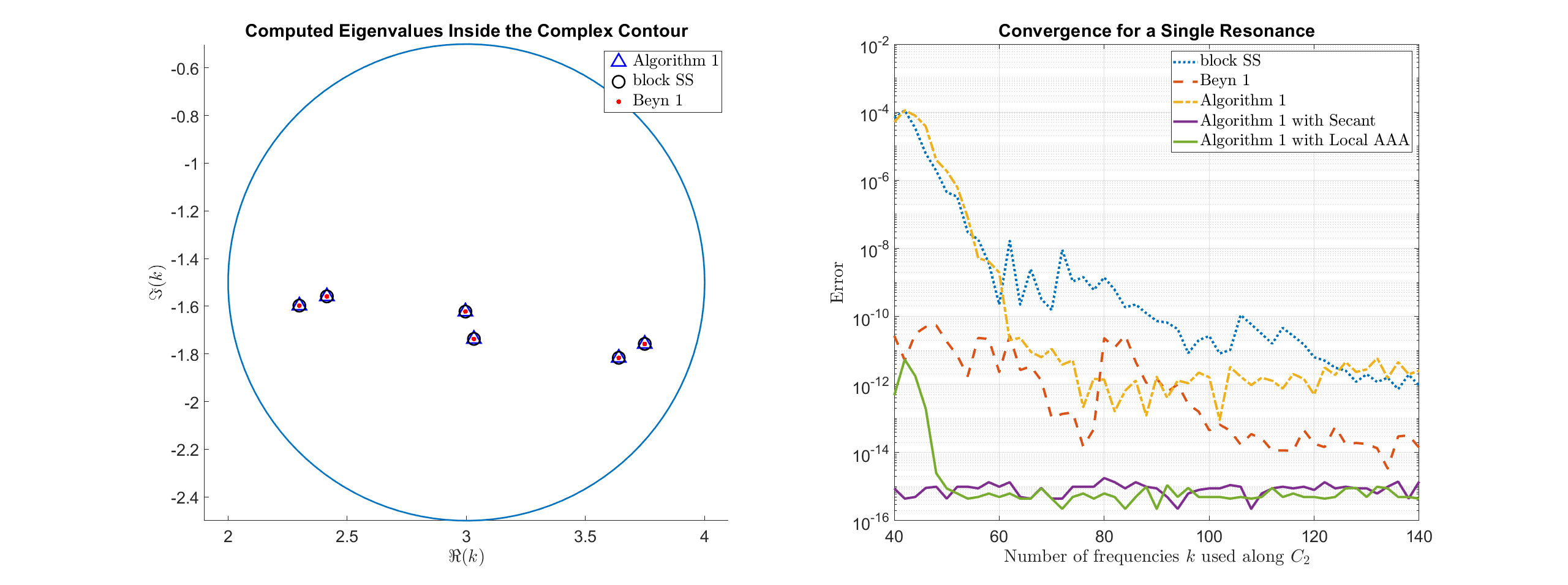}
        \caption{Comparison of Algorithm \ref{alg:one} to the block SS
          and Beyn 1 methods on the problem of evaluation of
          scattering poles outside the kite. Left panel: Eigenvalues
          computed by all three methods. Right panel: Errors 
          computed by comparison with a secant method evaluation of the
          eigenvalue near $2.299- 1.597i$. The curves labeled ``with
          secant'' were obtained by following the initial eigenvalue
          determination by four iterations of the secant method. For
          the curve labeled ``with local AAA'', four points were
          sampled on a circle of radius $1e-5$ around the initial AAA
          approximation of the pole together with a degree $1$
          rational approximant. To avoid underflow the maximum between
          the error and machine precision is plotted in all cases.}
    \label{fig:contourcomparison}
\end{figure} 
For our comparisons we consider the problem of evaluation of
eigenvalues of the Dirichlet problem in the exterior of the
kite-shaped curve used for the experiments in
Figure~\ref{fig:kites}. In detail, using Algorithm~\ref{alg:one} as
well as the block SS and Beyn~1 methods we evaluate the eigenvalues
contained within the circular contour $\mathcal{C}$ with center at
$3-1.5i$ and radius $1$, on the basis of equally spaced points along
$\mathcal{C}$. The left panel in Figure~\ref{fig:contourcomparison}
shows the eigenvalues within the curve $\mathcal{C}$ as computed by
each method, while the corresponding right panel displays the error
resulting from use of the various algorithms in the evaluation of the
eigenvalue near $2.299- 1.597i$ as a function of the number of points
used along the contour $\mathcal{C}$, using the computation by Algorithm~\ref{alg:one}
with secant method as the reference solution. Errors resulting from
use of the two different accuracy refinement methods introduced in
Section~\ref{sec:NEP}, namely, the secant method and the localized AAA
approximation, are included in the figure as well.

 \begin{figure}[h]
\begin{center}
\vskip 20pt
\includegraphics[scale=.7]{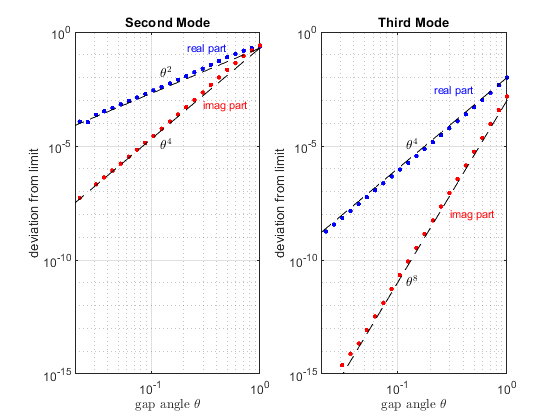}
\end{center}
\caption{Convergence of $\Re(k)$ and $\Im(k)$ to their limiting values
  as the gap size $\theta$ shrinks to zero for the nearly degenerate
  open circle eigenfunctions displayed as the second and third images
  on the top row of Figure~\ref{fig:circles}.  For both modes the
  imaginary parts (decay rates) converge at twice the rates of the
  real parts (spatial frequencies), and both the rates for the real
  and imaginary parts for the third mode are twice those for the
  second mode. Note that for the third mode the approximate nodal line
  is aligned with the gap and thus results in a weaker coupling of
  interior and exterior fields and associated faster
  convergence.}\label{gap-shrink}
\end{figure}

Figure~\ref{fig:contourcomparison} shows that Algorithm~\ref{alg:one}
and the block SS method converge at similar rates, at least for a
range of frequencies.  While Beyn~1 converges faster than
Algorithm~\ref{alg:one} and block SS without local refinement, use of
a local refinement technique on any of the three methods leads to
similar convergence rates for the eigenvalues.  Further, in situations
where the LU factorization of the matrix discretization of the
boundary integral equations is not available, such as may be the case
in problems requiring large surface- or even
volumetric-discretizations in 3D space, for which iterative solvers
may need to be used, Algorithm~\ref{alg:one} may prove more
advantageous than Beyn~1---for which the need to incorporate a total
of $L$ matrix solves per frequency $k$ along the contour, with values
of $L$ of the order of, say, a few tens, may become prohibitively
expensive.

As indicated in Section~\ref{intro}, the fact that
Algorithm~\ref{alg:one} does not utilize a quadrature rule gives rise
to a number of advantages in comparison with contour integration-based
methods. In particular, the use of rational approximation results in
geometric flexibility, such as, e.g., enabling the use of complex
frequencies $k$ that lie on a rectangular boundary in the complex plane,
or even simple segments along the real axis, which in turn allows for
simple adaptive algorithms such as the one introduced in
Section~\ref{sec:NEP} (Algorithm~\ref{alg:two}) without suffering
a deterioration in convergence. While the use of such rectangular
domains and associated adaptive algorithms could be envisioned in the
context of contour-integration based methods, such a strategy may
result in a significantly slower exponential rate of convergence than
a circular contour would and, if the trapezoidal rule is used over the
rectangular contour, where the contour integrand is not smooth, a
somewhat erratic convergence may result, as illustrated in~\cite[Fig.\
5.1(b)]{guttel2017nonlinear}. If Gauss-Lobatto quadrature is used
on each side of a rectangle, in turn, the exponent that characterizes
the exponential convergence rate is the approximately half of the one
associated with the trapezoidal rule on a circular contour~\cite[Sec.\
$5.2$]{guttel2017nonlinear}.  As a significant additional advantage,
Algorithm~\ref{alg:two} allows for the reuse of from previously used
data points as the frequency mesh is refined with the goal of
obtaining added accuracy for added accuracy---a strategy which cannot
be used in the context of the Gauss-Lobatto quadrature. Finally, the
applicability of the proposed algorithms on real segments instead of
closed complex contours provides an important benefit e.g. for
applications concerning real eigenvalues.

\subsection{Dependence on Gap Size}\label{sec:gapsize}

With a tool in hand to compute scattering poles quickly and
accurately, certain mathematical aspects of their behavior may be
considered, such as, e.g., the convergence rate of a complex
scattering pole associated with a open arc equal to the difference
between a closed curve and a gap section, as the gap size shrinks to
zero. As an example we consider the second and third modes shown in
the top row of Figure~\ref{fig:circles}, and associated eigenvalues
(but for varying gap sizes $\theta$) which, for small $\theta$ are
approximately equal to the lowest double eigenvalue the closed unit
circle ($\theta =0$)---namely the first root $\approx 3.8317$ of the
Bessel function $J_1(x)$.  
Figure~\ref{gap-shrink} demonstrates the convergence of $\Re(k)$ (the
``spatial frequency'') and $\Im(k)$ (the ``decay rate'') to their limiting
values as $\theta$ shrinks to zero.  The combined set of
convergence computations for the two modes considered was completed in
approximately four seconds on a single core of a present-day laptop
computer.

The first image in Figure~\ref{gap-shrink} shows that for the second
mode in Figure~\ref{fig:circles}, whose convergence rate we presume to
correspond to generic gap-shrinking eigenvalue convergence behavior,
the frequency converges to its limiting value $\approx 3.8317$ at the
rate $O(\theta^2)$ and the decay rate to its limiting value 0 at the
rate $O(\theta^4)$.  For the third mode in the figure, in turn, the
rates double to $O(\theta^4)$ and $O(\theta^8)$.  This is related to
the alignment of the nodal line of the eigenfunction with the gap,
which reduces the field coupling between the interior and exterior
regions of the open cavity. We are not aware of any theoretical or
computational studies reporting on such convergence rates for
scattering poles as opening gaps tend to zero.

The computations of Figs.~\ref{fig:circles} and~\ref{gap-shrink} correspond to Dirichlet boundary
 conditions, and it is interesting to ask the same questions in the
Neumann case.  Here the first degenerate eigenvalue for the circle,
associated with the first zero of the derivative of $J_1(r)$,
 is $1.84118\dots$.  Once again, introducing a gap in the boundary
 breaks the degeneracy.  In separate computations (not shown)
we have found that for the mode whose nodal line is aligned
 with the gap, the eigenvalue now converges at the rates $\theta^2$
for the real part and $\theta^4$ for the imaginary part, whereas
for the more generic mode whose nodal line is orthogonal to the gap,
 the rates slow down dramatically to $1/|\log\theta|$ and $1/|\log\theta|^2$,
 respectively.
 
\subsection{High-frequency Examples}\label{sec:highfreq}
Several higher-frequency examples are considered in what follows. To
verify the accuracy of the proposed methods in high-frequency cases
and in regions containing large numbers of eigenvalues we applied the
real line version of algorithm $2$ to obtain all $1244$ distinct
interior Laplace eigenvalues of the unit circle that are contained in
the interval $[1,100]$. The algorithm automatically obtained all
$1244$ eigenvalues with near machine accuracy in an average
computational time of $1.09$ seconds per eigenvalue in a single CPU
core. This relatively slow figure is dominated by the higher
eigenvalues; for example the average time per eigenvalue in the
interval $[0,25]$ is only about $0.08$ seconds.

\begin{figure}[h]
    \centering
     \includegraphics[scale=0.5]{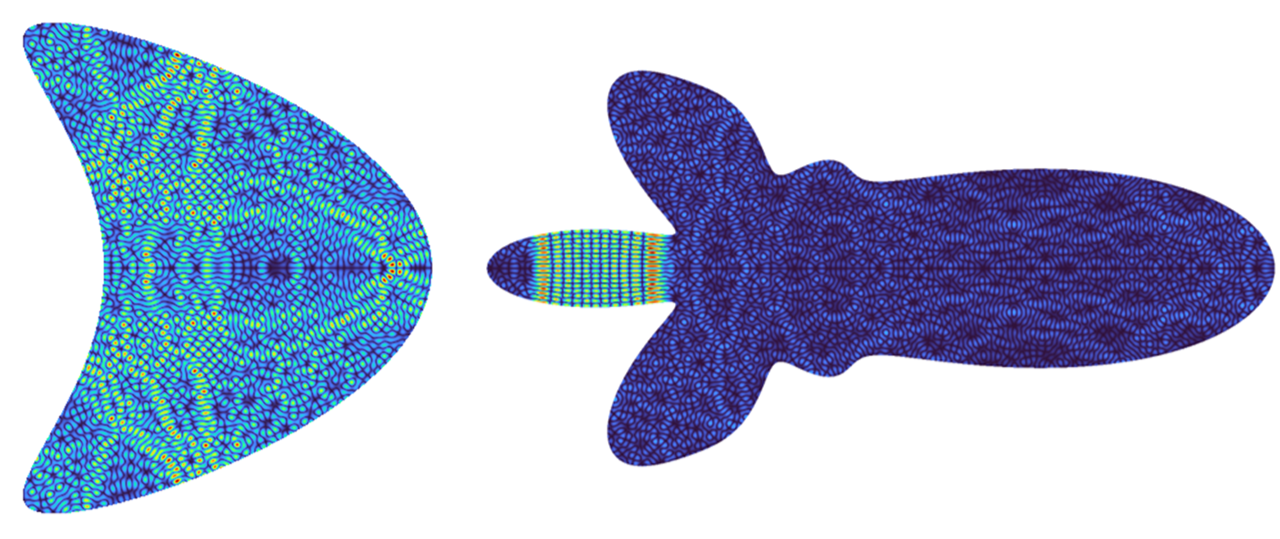}
     \caption{Left: Interior eigenfunction for the kite-shaped domain,
       with eigenvalue $100.1846738596$. Right: Interior eigenfunction
       for the rocket-shaped domain, with eigenvalue
       $399.9730212127$. All digits listed are believed correct on the
       basis of numerical convergence analyses.}
    \label{fig:Kite}
\end{figure}

\begin{figure}[p]
    \centering
    \includegraphics[width=\textwidth]{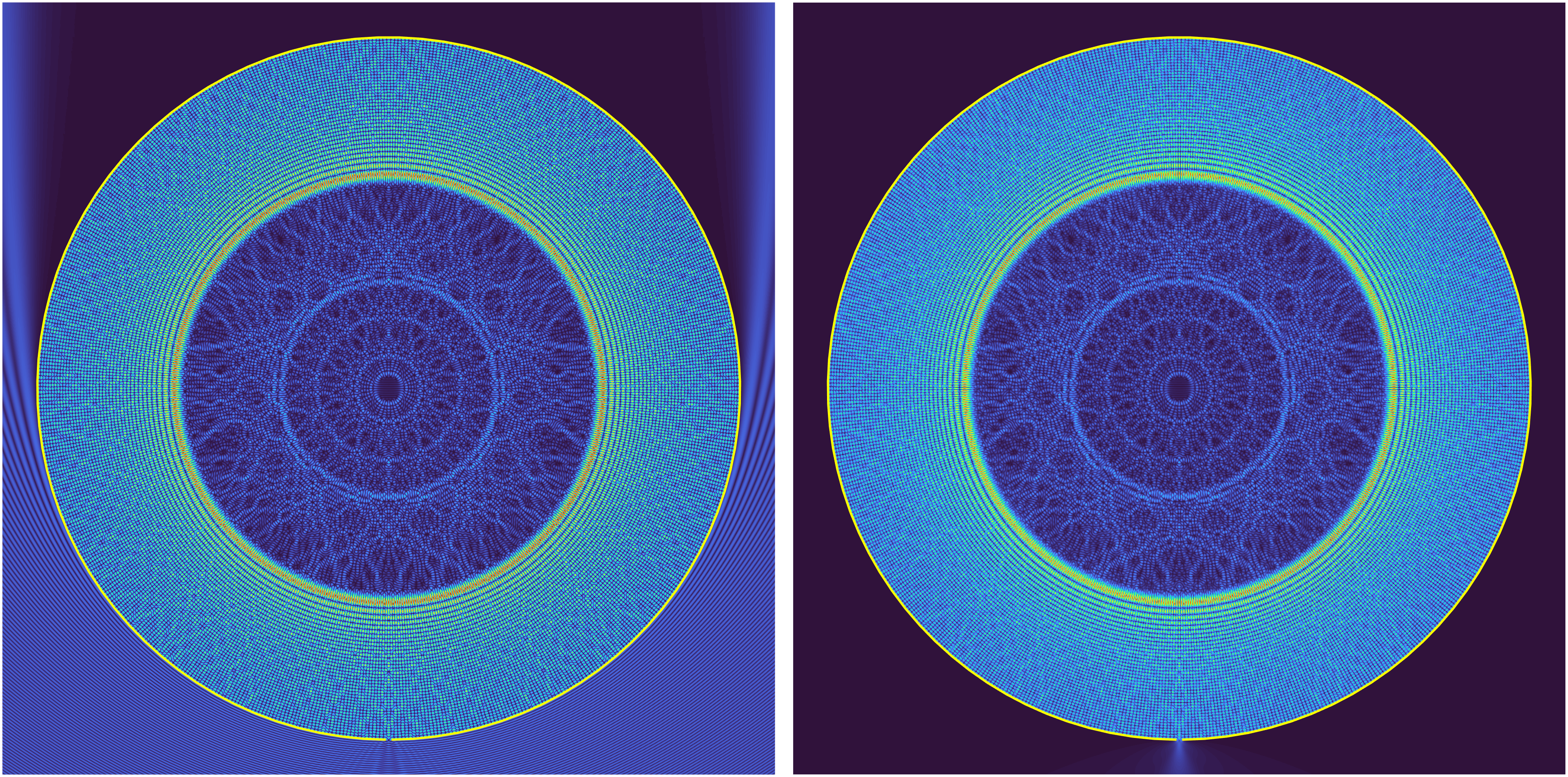}
    \caption{Comparison of the solution of the scattering problem with
      $k = 499.9073989141$ under vertical and upward plane-wave
      incidence (left) and the eigenfunction corresponding to the
      eigenvalue $ 499.9073989141-0.000779974959i$ (right) for an open
      circular cavity with aperture size $\pi/100$. All digits listed are
      believed correct on the basis of numerical convergence analyses.}
    \label{fig:Circle}
\end{figure}

\begin{figure}
    \centering
    \includegraphics[width=\textwidth]{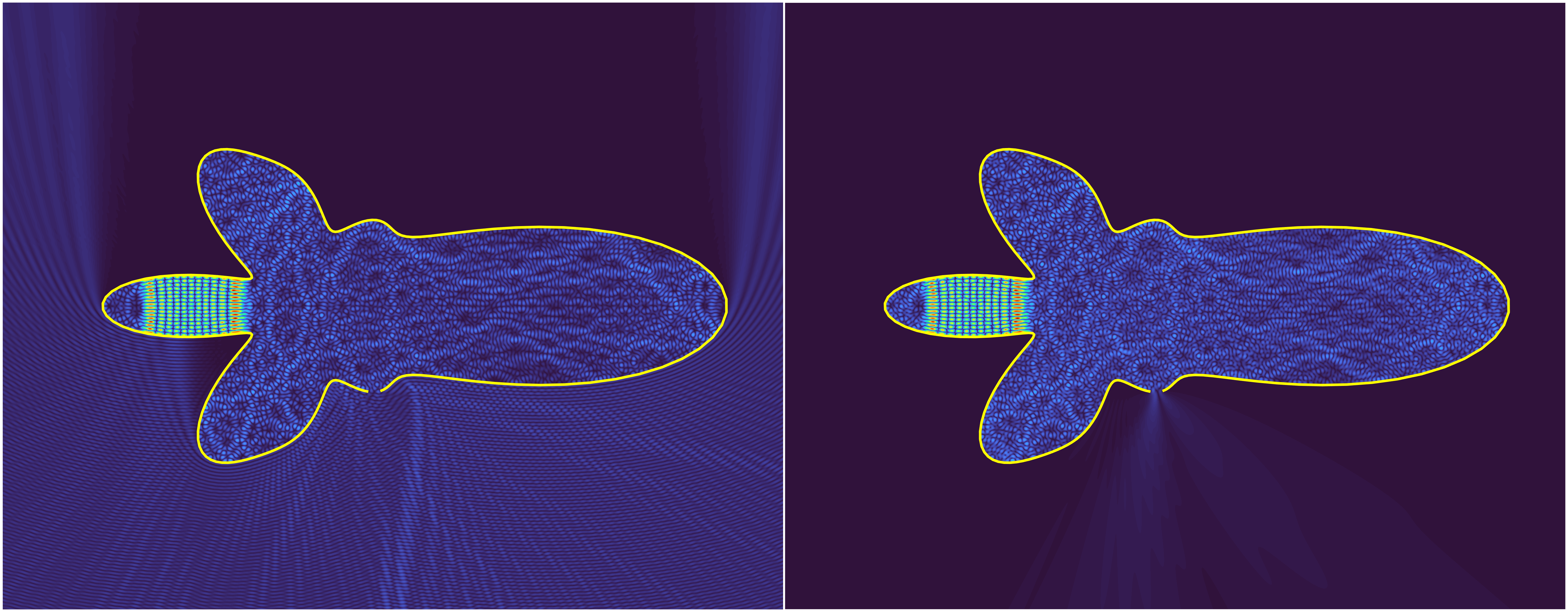}
    \caption{Comparison of the solution of the scattering problem with
      $k = 399.9694808817$ under vertical plane-wave incidence (left)
      and the eigenfunction corresponding to the eigenvalue
      $399.9694808817 - 0.00434495360i$ (right) for a rocket-like open
      cavity. The opening in the right curve is approximately equal to
      $0.6\%$ of the curve length. All digits listed are believed
      correct on the basis of convergence analyses.}
    \label{fig:Rockets}
\end{figure}

Having verified our algorithm in the high-frequency regime, we
subsequently applied the proposed methods to several high-frequency
eigenvalue problems for interior and open-cavity
domains. Figure~\ref{fig:Kite} displays high-frequency interior
eigenfunctions for kite-shaped and rocket-shaped cavities. The right
panels in Figures \ref{fig:Circle} and \ref{fig:Rockets}, in turn,
present eigenfunctions for circular and rocket-shaped open cavities,
while the left panels of these figures present the solutions of
problems of scattering by the same cavities under vertical
upward-facing plane-wave illumination at frequencies equal to the real
parts of the eigenvalues associated with the corresponding right
panels. Clearly, the left panel scattering solutions bear a close
resemblance with the eigenfuncions displayed in the right panels,
suggesting that the open-cavity eigenfunctions can be excited by
adequately oriented incident fields penetrating through the small
apertures.  All digits displayed in the figure captions have been
found to be correct by means of studies of convergence in both the
discretization of the integral equations used and of the eigensearch
algorithms utilized. For reference in Figures~\ref{fig:Circle} and
\ref{fig:Rockets}, the most refined solution used $16$ points per
wavelength which led to matrix discretizations of size
$5000\times 5000$ and $4000 \times 4000$ respectively and produced
13-digit accuracies.

\section{Conclusions}\label{sec:conclusion}
This paper introduced a novel numerical algorithm for the evaluation
of real and complex eigenvalues and eigenfunctions of general NEPs,
including NEPs associated with Laplace eigenvalue and
scattering-resonance problems for open and closed domains. Based on
use of adaptive eigenvalue search methods based on AAA rational
approximation combined with secant-method refinement, the algorithm
produces highly accurate eigenpairs for challenging eigenproblems at
both low and high frequencies.  Comparisons with well-known
contour-integration methods demonstrated a number of advantages of the
proposed approach, including fast convergence and geometric
flexibility. The latter characteristic is significant, in that it
greatly facilitates use of the algorithm in a rectangular
refinement-based adaptive strategy---resulting in automatic evaluation
of all eigenvalues contained within a given region with near machine
precision accuracy.

\section*{Acknowledgements}
OB gratefully acknowledges support from NSF and AFOSR under contracts
DMS-2109831 and FA9550-21-1-0373. MS acknowledges support from the
National Science Foundation Graduate Research Fellowship under Grant
No.\ 2139433. The authors gladly acknowledge helpful discussions with
the following individuals: Stefan G\"uttel, Jonathan Keating, Yuji
Nakatsukasa, Ory Schnitzer, Euan Spence, and Maciej Zworski.

\bibliographystyle{plain}
\bibliography{main}
\end{document}